\documentclass[10pt, onecolumn]{article}  % Force single-column layout

\usepackage{graphicx}       % For images
\usepackage{hyperref}       % For metadata and links
\usepackage{url}            % For proper URL formatting
\usepackage{amsmath}        % For math equations
\usepackage{amsthm}
\usepackage{amssymb}        % For math symbols
\usepackage{booktabs}       % For professional tables
\usepackage{xcolor}         % For color control
\usepackage{lipsum}         % For placeholder text (remove in actual document)
\usepackage{tabularx}
\usepackage{geometry}
\usepackage{wrapfig}
\usepackage{fancybox,fancyhdr}
\usepackage{lastpage}
\usepackage{multirow}
\usepackage{amscd}
\usepackage{indentfirst}

\title{\bfseries\Large On an analogue of Frobenius formalism for 3-algebras and pentagon equations solutions arising from projectors}
\author{
  Ramil K. Aliev\\
  \\
  Lomonosov MSU, Leninskie Gory, 1, 119991, Moscow, Russia  \\
  \texttt{alievrk@yandex.ru}
}
\date{}  % Omit date to avoid rebuild issues

\hypersetup{
  colorlinks=true,
  linkcolor=blue,
  citecolor=red,
  urlcolor=magenta,
  pdfinfo={
    Title={On an analogue of Frobenius formalism for 3-algebras and pentagon equations solutions arising from projectors},
    Author={Ramil K. Aliev},
    Keywords={3-algebra, Frobenius compatibility, 3-manifold, ternary, pentagon equations, $F$-matrix},
    CreationDate={},
    ModDate={}
  }
}

\newcommand{\keywords}[1]{\par\noindent\textbf{Keywords:} #1}

\begin{document}
% ===== Title Section =====
\maketitle

% ===== Abstract & Keywords =====
\begin{abstract}
  Ruth J.Lawrence introduced a notion of a 3-algebra to construct invariants of 3-manifolds based on their triangulations in her paper "Algebras and triangular relations". Her primary definition is suitable for certain triangulations only although a hint to handle arbitrary ones has been proposed. Here I introduce an analogue of Frobenius compatibility for one class of 3-algebras hence I obtain a way to construct a full 3-algebra. Additionally, I provide with examples of a 3-algebra and invariants for lens spaces. Moreover, it leads to a new family of pentagon equations solutions: arising from projectors.
  
\end{abstract}
\keywords{3-algebra, Frobenius compatibility, 3-manifold, ternary, pentagon equations, $F$-matrix}  % Add 3-5 relevant keywords

% ===== Main Content =====
\section{Introduction}

\subsection{3-algebra}

\textbf{Definition. } A \textbf{3-algebra}, $A$, over a field $K$, is a vector space over $K$ with a basis $\{e_i\}_{i = 1}^{n}$, endowed with $K$-linear maps:

$$ P: A  \rightarrow A \quad (P^3 = id),$$

$$ m: A^{\otimes 3} \rightarrow  A,$$

$$ \overline{m}: A^{\otimes 2} \rightarrow  A^{\otimes 2},$$

satisfying the following axioms:

(i) $m(m \otimes 1 \otimes 1) = m(1 \otimes 1 \otimes m)\sigma_{34}(1 \otimes \overline{m} \otimes 1 \otimes 1)\sigma_{34}$,

(ii) $(1 \otimes m)\sigma_{23}(\overline{m} \otimes 1 \otimes 1)=\overline{m}(1 \otimes m)\sigma_{12}(P^{-1} \otimes 1 \otimes 1 \otimes 1)(\overline{m} \otimes 1 \otimes 1)(P \otimes P \otimes 1 \otimes 1)\sigma_{23}$,

(iii) $\overline{m}(m \otimes 1)=(1 \otimes m)\sigma_{12}(P^{2} \otimes \overline{m} \otimes 1)(1 \otimes 1 \otimes \overline{m})\sigma_{12}\sigma_{23}$,

(iv) $(1 \otimes \overline{m})\sigma_{12}(1 \otimes \overline{m})=(\overline{m} \otimes 1)(1 \otimes \overline{m})(P \otimes P \otimes 1)(\overline{m} \otimes 1)(1 \otimes P^{-1} \otimes 1)$,

(v) $(1 \otimes m)\sigma_{23}(\overline{m} \otimes P^{2} \otimes 1)=(m \otimes 1)(1 \otimes 1 \otimes \overline{m})$,

(vi) $Pm=m(P \otimes P \otimes P)\sigma_{23}\sigma_{12}$,

(vii) $\overline{m}(P^{2} \otimes P)\sigma_{12}=\sigma_{12}\overline{m}(P^{2} \otimes P)$

($\sigma_{ij}$ is a map permuting $e_i$ and $e_j$).

$e_{i_1 ... i_n}$ stands for $e_{i_1} \otimes ... \otimes e_{i_n}$.

Operations $m$ and $\overline{m}$ are visualised by tetrahedra: each face of a tetrahedron is assigned "-" (for inputs) or "+" (for outputs). A copy of $A$ or $A^*$ is placed on each face, for "-" and "+" respectively.

To better understand, see figures 1-7. They represent the "view above" on input faces.

\begin{figure}[!ht]
\begin{center}
\includegraphics[scale=0.9]{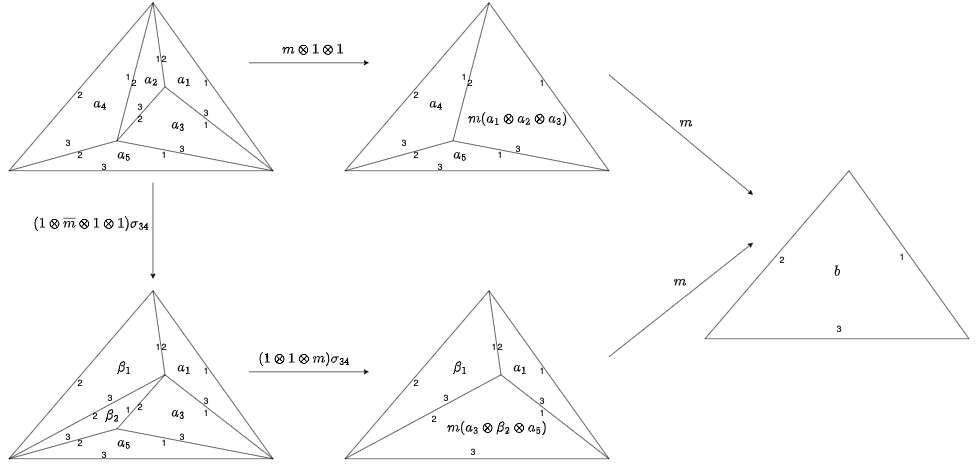}\caption{Axiom (i)}\label{figure1}
\end{center}
\end{figure}

\clearpage

\begin{figure}[!ht]
\begin{center}
\includegraphics[scale=0.8]{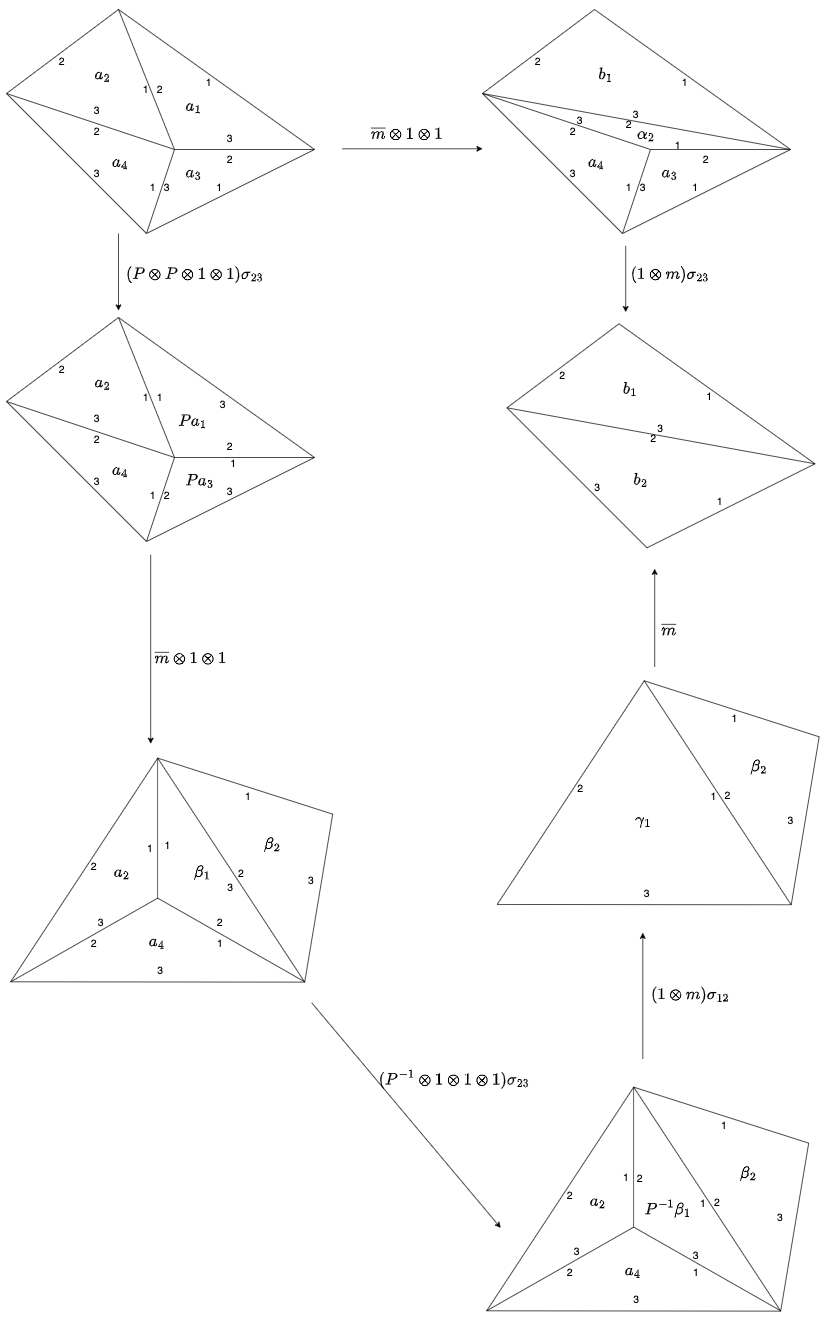}\caption{Axiom (ii)}\label{figure2}
\end{center}
\end{figure}

%%\clearpage

\begin{figure}[!ht]
\begin{center}
\includegraphics[scale=0.5]{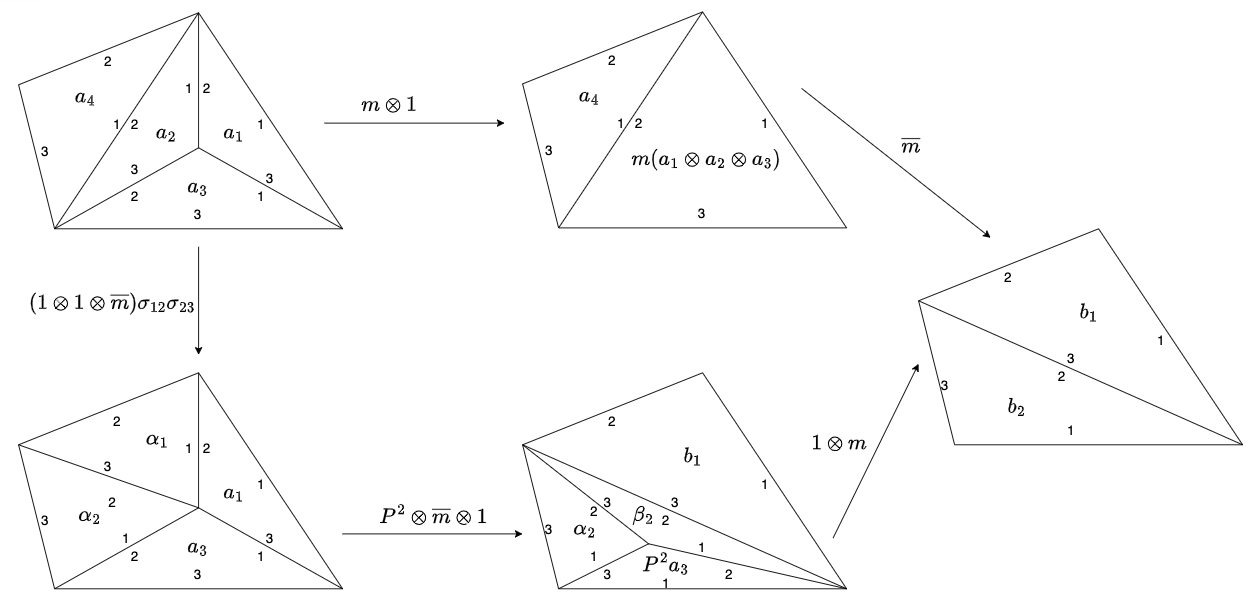}\caption{Axiom (iii)}\label{figure3}
\end{center}
\end{figure}

%%\clearpage

\begin{figure}[!ht]
\begin{center}
\includegraphics[scale=0.5]{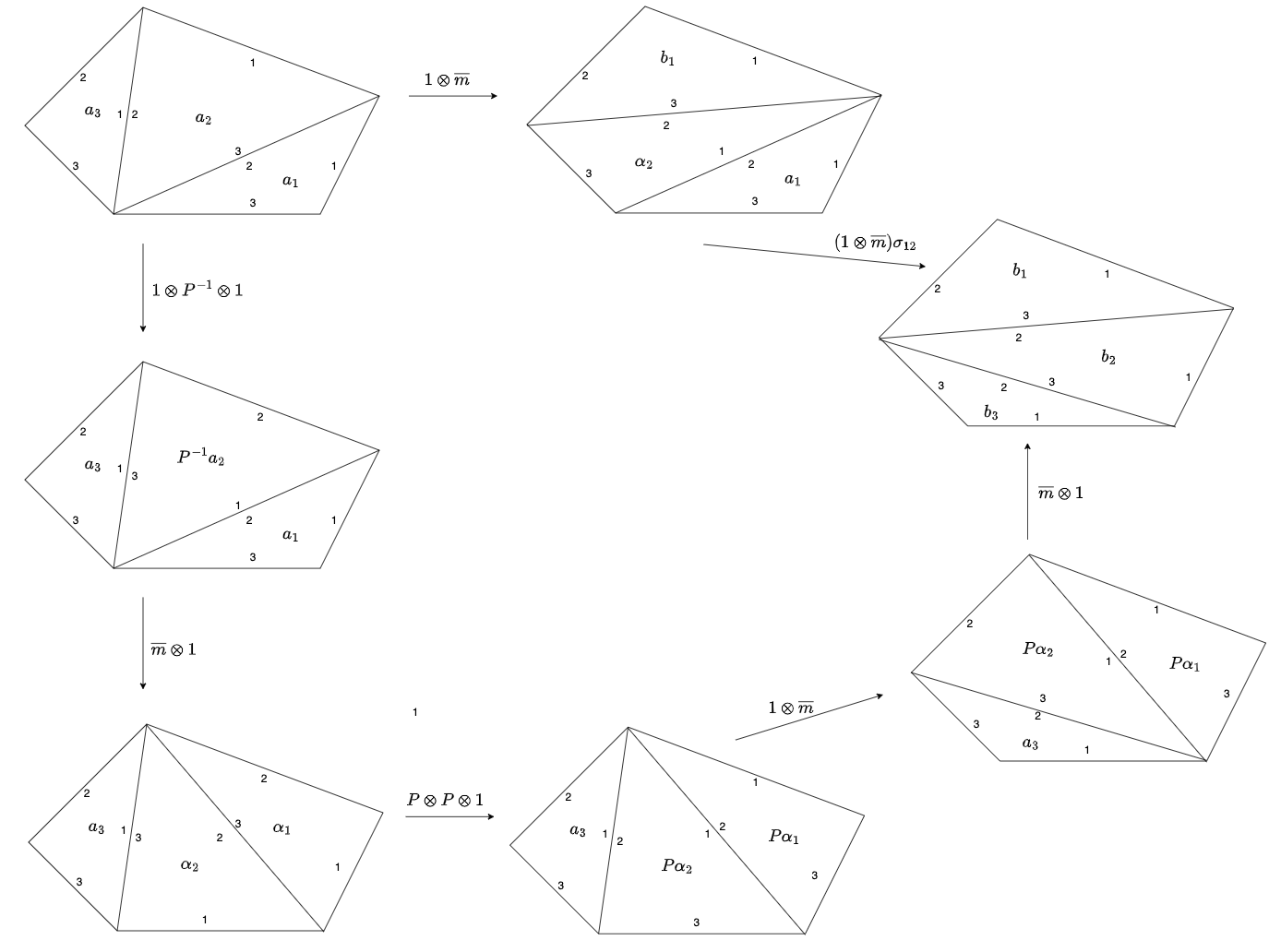}\caption{Axiom (iv)}\label{figure4}
\end{center}
\end{figure}

%%\clearpage

\begin{figure}[!ht]
\begin{center}
\includegraphics[scale=0.6]{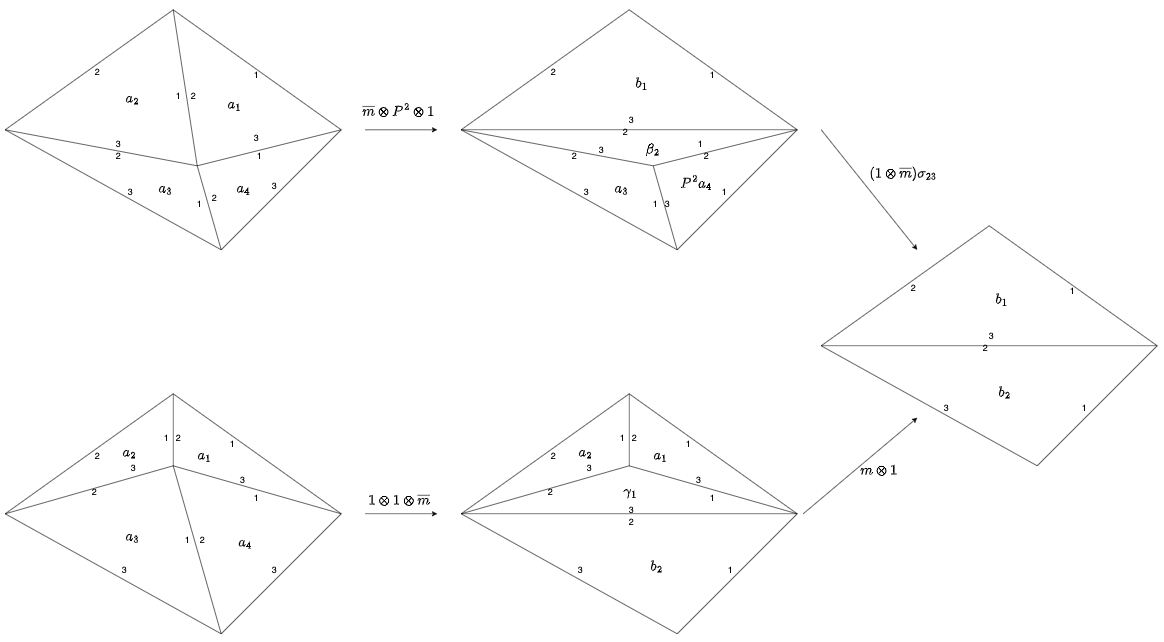}\caption{Axiom (v)}\label{figure5}
\end{center}
\end{figure}

%%\clearpage

\begin{figure}[!ht]
\begin{center}
\includegraphics[scale=0.9]{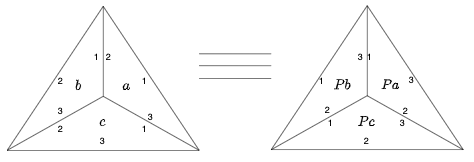}\caption{Axiom (vi)}\label{figure6}
\end{center}
\end{figure}

\begin{figure}[!ht]
\begin{center}
\includegraphics[scale=0.9]{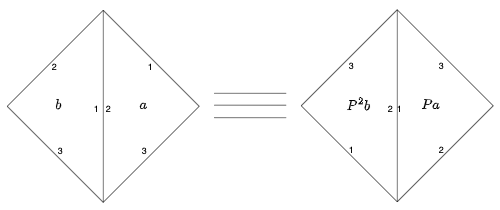}\caption{Axiom (vii)}\label{figure7}
\end{center}
\end{figure}
 
\clearpage
 
Let us introduce structural constants for $\overline{m}$ and $m$:

$$\overline{m}(e_i \otimes e_j) = Q_{ij}^{st}e_{st},$$
$$m(e_i \otimes e_j \otimes e_k) = Q_{ijk}^{t}e_{t}.$$

\subsubsection{Map $P$}

Inside each triangular face, edges are labelled 1,2 or 3. An edge may have different labels from incident faces. If a face has 1-2-3 on its edges in counterclockwise order, it is assigned "-", in clockwise -- "+". Here, (counter)clockwise order of labels on a face is meant the one seen "from the outside of a tetrahedron". 

Labels on edges denote the order of arguments corresponding to adjacent faces.

The map $P$ acts on 3-algebra $A$. It is of order 3. Its role is to rotate a triangular face to get certain labels distribution so that this face matches the order of maps $\overline{m}$ and ${m}$ described later.

!These labels shouldn't be confused with "naming by vertices" fixed once and for all!

\begin{figure}[!ht]
\begin{center}
\includegraphics[scale=0.8]{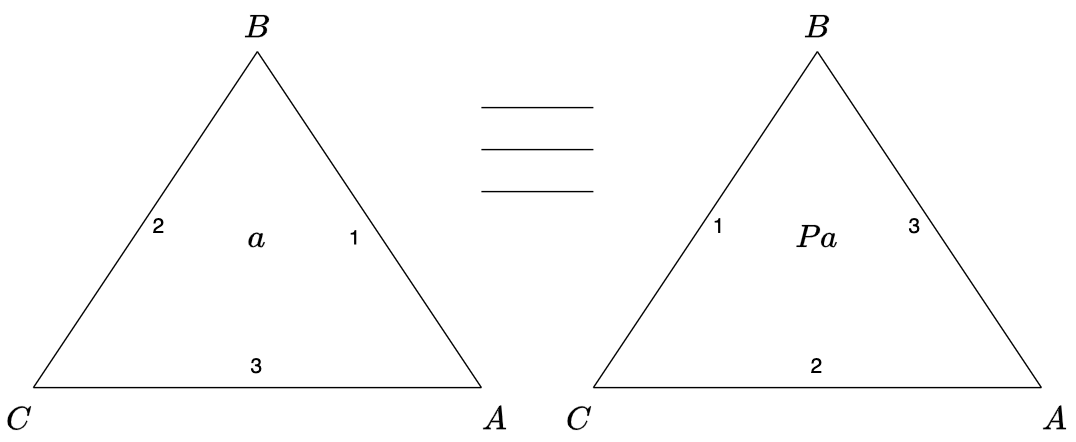}\caption{Map $P$}\label{figure8}
\end{center}
\end{figure}

\subsubsection{Map $m$}

A tetrahedron associated with map $m$ has three input faces and one output face. $i-$th argument for $m$ corresponds to the face the output face shares an edge with label $i$ with.

\begin{figure}[!ht]
\begin{center}
\includegraphics[scale=0.8]{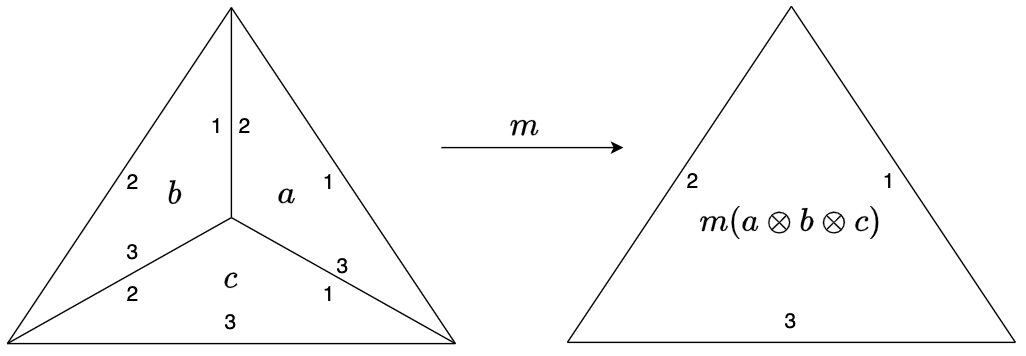}\caption{Map $m$ (The LHS picture is the view from the outside of a tetrahedron, the RHS -- from the inside)}\label{figure9}
\end{center}
\end{figure}
 
\subsubsection{Map $\overline{m}$}

A tetrahedron associated with map $\overline{m}$ has two input faces and two output faces. 

Its input and output faces are ordered as follows:

1. Take two input faces. 

2. Look at the edges they don't share.

3. There're only one with label 1 and only one with label 2. Remember them.

4. The face incident to the edge with label 1 is the first input, to the edge with label 2 -- the second input.

5. The output face containing this edges is the first output.

\begin{figure}[!ht]
\begin{center}
\includegraphics[scale=0.7]{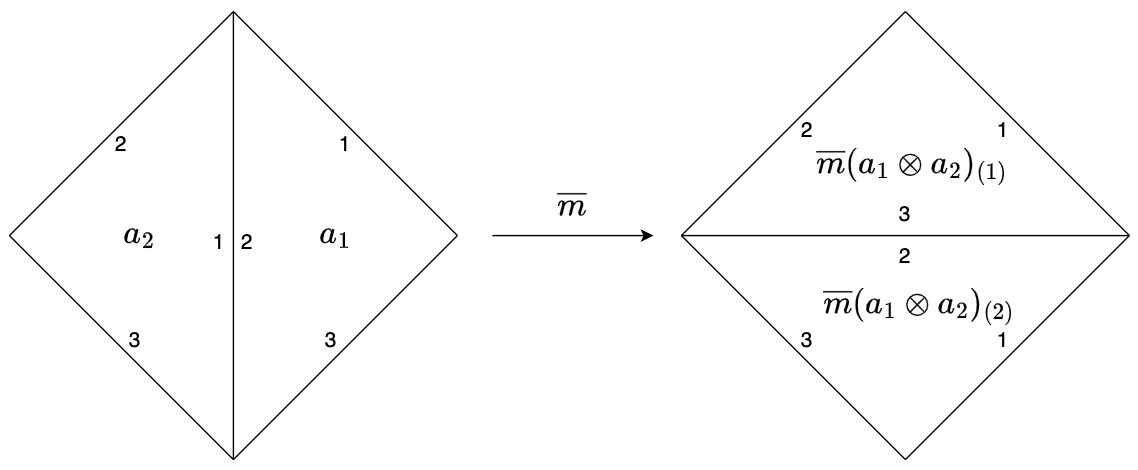}\caption{Map $\overline{m}$  (The LHS picture is the view from the outside of a tetrahedron, the RHS -- from the inside)}\label{figure10}
\end{center}
\end{figure}

\subsubsection{Construction of the invariant}

A 3-manifold is split into tetrahedra, each tetrahedron is assigned an operation, $m$ or $\overline{m}$ depending on how its faces are labelled (here is meant the original notion of a 3-algebra which allows three-input-one-output and two-input-two-output tetrahedra only). A triangular face incident to two tetrahedra is input for one of them and output for another. Such a face denotes contraction of its tetrahedra's operations by its corresponding argument. This intermediate pair of tetrahedra is assigned a composition of operations assigned to these tetrahedra.

If the manifold has a boundary, each of its component is considered input or output. In its turn, its shelling consists of triangular faces standing for $A$ or $A^*$. These are the inputs and outputs for the operation corresponding to the manifold.

If the manifold is closed, one splits it and obtains submanifolds with boundary, deals with them the way stated above, then glues these submanifolds and contracts respective operations by arguments corresponding to common faces the manifold was partitioned along.

Axioms (i)-(iv) formalise 2-3 Pachner move, axiom (v) ensures that the result is independent on internal vertices, relations (vi) and (vii) axiomatise symmetries of operations $m$ and $\overline{m}$ respectively. 

\subsection{Frobenius compatibility}

\textbf{Definition. } A multiplication in an algebra $A$ and a metric $(\cdot, \cdot)$ are \textbf{(Frobenius-)compatible} if $\forall a,b,c \in A$  $(ab,c)=(a,bc)$.

\textbf{Definition. } {Frobenius algebra} is an associative unital algebra endowed with metric compatible with multiplication.

\section{Frobenius 3-algebra and full 3-algebra}

As R.Lawrence mentioned, a 3-algebra may be extended to a so called \textbf{full} 3-algebra in order to allow arbitrary triangulations. A full 3-algebra is supposed to have all maps $m_j : A^{\otimes j} \rightarrow  A^{\otimes (4-j)}, \quad (j=0,1,2,3)$. Nevertheless, no explicit definition or method to construct one has been suggested.

If signs on internal faces don't affect the result, one should have an opportunity to change a sign with respect to the same tetrahedron, i.e. transforming an output to input and vice-versa. Hence the idea to generalise Frobenius compatibility to 3-algebras arises.

However, only a certain class of 3-algebras is considered here: if $P$ is a permutation of basis vectors. Let $e_1, ..., e_n$ be a basis of 3-algebra. For conveniency, let us write $P(e_i)$ as $e_{P(i)}$. For these 3-algebras, somewhat analogue of Frobenius compatibility may be defined.

In this case, axioms (vi) and (vii) look like the following:

(vi) $Q_{ijk}^s = Q_{P(j)P(k)P(i)}^{P(s)}$,

(vii) $Q_{ij}^{st} = Q_{P^2(j)P(i)}^{P^2(t)P(s)}$.

\textbf{Definition. } Operations $\overline{m}$ and $m$ are \textbf{compatible with respect to bilinear form $h(\cdot, \cdot)$} if the following holds:

$$Q_{ijk}^s = Q_{ij}^{st}h_{tk} \quad \forall i,j,k,s,t. \quad (1)$$

\textbf{Proposition. }  Let $A$ be a 3-algebra with $P$ being a permutation on basis vectors. If an non-degenerate bilinear form $h(\cdot, \cdot)$ satisfies the condition

$$h_{jk} = h_{P(k)P^2(j)} \quad \forall j,k, \quad (2)$$

the compatibility between $\overline{m}$ and $m$ with respect to this form results in the following conditions:

a) the axiom (v) holds

b) axioms (i)-(iv) are equivalent

\begin{proof}

a) Axiom (v), applied to $e_{ijkl}$:

$$Q_{ijk'}^s Q_{kl}^{k't} = Q_{ij}^{sj'}Q_{P^2(k)j'l}^{t}.  \quad (3)$$

Insert a bilinear form $h$ into the LHS of (3) and consider the compatibility (1):

$$Q_{ijk'}^s (h^{k'j'}h_{j'k'})Q_{kl}^{k't} = Q_{ij}^{sj'}(h_{j'k'}Q_{kl}^{k't}).  \quad (4)$$

In the RHS of (4) swap indices in the expression inside the parentheses using axiom (vii) and property (2):

$$h_{j'k'}Q_{kl}^{k't} = h_{P(k')P^2(j')}Q_{P^2(l)P(k)}^{P^2(t)P(k')} = Q_{P^2(l)P(k)}^{P^2(t)P(k')}h_{P(k')P^2(j')} = Q_{P^2(l)P(k)P^2(j')}^{P^2(t)}.$$

Finally, apply axiom (vi):

$$Q_{P^2(l)P(k)P^2(j')}^{P^2(t)} = Q_{P^2(k)j'l}^{t}.$$

Hence, axiom (v) holds.

b) Expand axioms (i)-(iv) in coordinate form:

1) Axiom (i), applied to $e_{ijkst}$:

$$Q_{l'st}^{r}Q_{ijk}^{l'} = Q_{ij'l}^{r}Q_{ks't}^{l}Q_{js}^{j's'}.  \quad (5)$$

In the RHS of (5), consider axiom (v) and replace $Q_{ks't}^{l}Q_{js}^{j's'}$ with $Q_{jsk'}^{j'}Q_{P(k)t}^{k'l}$:

$$Q_{l'st}^{r}Q_{ijk}^{l'} = Q_{ij'l}^{r}Q_{jsk'}^{j'}Q_{P(k)t}^{k'l}.  \quad (6)$$

Using axiom (vi), replace the constants of $m$ with equal ones and substitute them into (6):

$Q_{l'st}^{r} = Q_{P^2(t)P^2(l')P^2(s)}^{P^2(r)}, Q_{jsk'}^{j'}  = Q_{P^2(k')P^2(j)P^2(s)}^{P^2(j')}, Q_{ijk}^{l'} = Q_{P(j)P(k)P(i)}^{P(l')}, Q_{ij'l}^{r}  = Q_{P(j')P(l)p(i)}^{P(r)},$

$$ Q_{P^2(t)P^2(l')P^2(s)}^{P^2(r)}Q_{P(j)P(k)P(i)}^{P(l')} = Q_{P(j')P(l)P(i)}^{P(r)}Q_{P^2(k')P^2(j)P^2(s)}^{P^2(j')}Q_{P(k)t}^{k'l}.  \quad (7)$$

Taking the scalar product of (7) with $h^{P(i)q}$ and $h^{P^2(s)u}$, raise indices $P(i)$ and $P^2(s)$ in both sides of (7):

$$ Q_{P^2(t)P^2(l')}^{P^2(r)u}Q_{P(j)P(k)}^{P(l')q} = Q_{P(j')P(l)}^{P(r)q}Q_{P^2(k')P^2(j)}^{P^2(j')u}Q_{P(k)t}^{k'l}.  \quad (8)$$

2) Axiom (ii), applied to $e_{ijkl}$:

$$ Q_{ij}^{tj'}Q_{kj'l}^{s} = Q_{k's'}^{ts}Q_{P(i)P(k)}^{i'k'}Q_{P^2(i')jl}^{s'}.  \quad (9)$$

Taking the scalar product of (9) with $h^{lq}$, raise index $l$ in both sides of (9):

$$ Q_{ij}^{tj'}Q_{kj'}^{sq} = Q_{k's'}^{ts}Q_{P(i)P(k)}^{i'k'}Q_{P^2(i')j}^{s'q}.  \quad (10)$$

3) Axiom (iii), applied to $e_{ijkl}$:

$$ Q_{s'l}^{qr}Q_{ijk}^{s'} = Q_{P^2(k)t'l'}^{r}Q_{ij'}^{qt'}Q_{jl}^{j'l'}.  \quad (11)$$

Using axiom (vi), make a replacement and substitute into (11):

$Q_{P^2(k)t'l'}^{r}= Q_{P(t')P(l')k}^{P(r)},$

$$ Q_{s'l}^{qr}Q_{ijk}^{s'} = Q_{P(t')P(l')k}^{P(r)}Q_{ij'}^{qt'}Q_{jl}^{j'l'}.  \quad (12)$$

Taking the scalar product of (12) and $h^{kt}$, raise index $k$ in both sides of (12):

$$ Q_{s'l}^{qr}Q_{ij}^{s't} = Q_{P(t')P(l')}^{P(r)t}Q_{ij'}^{qt'}Q_{jl}^{j'l'}.  \quad (13)$$

4) Axiom (iv), applied to $e_{ijk}$;

$$ Q_{ik'}^{st}Q_{jk}^{rk'} = Q_{P(i')q'}^{rs}Q_{P(j')k}^{q't}Q_{iP^2(j)}^{i'j'}.  \quad (14)$$

Relations (8),(10),(13),(14), corresponding to axioms (i)-(iv), coincide up to renaming the indices thus the proposition is proven.

\end{proof}

Hence, axioms (i)-(iv) turn out to be equivalent not only in formal-philosophic sense (different 2-3 Pachner move formalisation, [3]) but also mathematically strictly.

Previous observations allows to introduce a notion of a Frobenius 3-algebra.

\textbf{Definition. } A tetraplet $(A,P,\overline{m},h)$ where:

-- $A$ -- a vector space over a field $K$

-- $P$ -- a permutation on basis vectors of $A$, $P^3 = id$

-- $\overline{m}:  A^{\otimes^2}  \rightarrow  A^{\otimes^2} $ -- a $K$-linear operator on  $A \otimes A$

-- $h$ -- non-degenerate bilinear form on $A$, satisfying the property:

$$h_{jk} = h_{P(k)P^2(j)} \quad \forall j,k,$$

is called a \textbf{Frobenius 3-algebra} if it fulfills the following relations:

(i$^*$) $(1 \otimes \overline{m})\sigma_{12}(1 \otimes \overline{m})=(\overline{m} \otimes 1)(1 \otimes \overline{m})(P \otimes P \otimes 1)(\overline{m} \otimes 1)(1 \otimes P^{-1} \otimes 1),$

(ii$^*$) $\overline{m}(P^{2} \otimes P)\sigma_{12}=\sigma_{12}\overline{m}(P^{2} \otimes P).$

A full 3-algebra is supposed to have all operations $m_j : A^{\otimes j} \rightarrow  A^{\otimes (4-j)}, \quad (j=0,1,2,3)$. Now let us consider the explicit definition of a full 3-algebra based on the notion of a Frobenius 3-algebra.

\textbf{Definition. } A Frobenius 3-algebra $(A,P,\overline{m},h)$ is called \textbf{full} if it has all maps $m_j : A^{\otimes j} \rightarrow  A^{\otimes (4-j)}, \quad (j=0,1,2,3)$ where:

-- $m_{2,2} = \overline{m}, \quad \overline{m}(e_{ij}) = Q_{ij}^{st}e_{st},$

-- $m_{3,1}(e_{ijk})= Q_{ijk}^{s}e_s, \quad Q_{ijk}^t = Q_{ij}^{st}h_{tk},$

-- $m_{1,3}(e_i) = Q_i^{stu}e_{stu}, \quad Q_i^{stu}:=Q_{ij}^{st}h^{ju},$

-- $m_{0,4}(1) = Q^{ijkl}e_{ijkl}, \quad Q^{ijkl}:=Q_{s}^{ijk}h^{sl},$

-- $m_{4,0}(e_{ijkl}) = Q_{ijkl}, \quad Q_{ijkl}:=Q_{ijk}^{t}h_{tl}.$

( In these terms, $m_{3,1} = m$ from Lawrence's definition)

Since $m_{3,1}$ and $m_{2,2}$ are connected via a non-degenerate bilinear form $h$, one can use it to define the missing operations $m_{j,4-j}: A^{\otimes j} \rightarrow  A^{\otimes (4-j)} \quad (j=0,1,4)$.

%Extra operations $m_{1,3}$, $m_{0,4}$ and $m_{4,0}$ fit axioms of a 3-algebra or a Frobenius 3-algebra.

%\textbf{Proposition. } Extra operations $m_{1,3}$, $m_{0,4}$ and $m_{4,0}$ fit axioms of a 3-algebra or a Frobenius 3-algebra.

%\begin{proof}

%Take the scalar product of $h_{ru}$ with (5) expressing applied to $e_{ijkst}$ axiom (i):

%$$Q_{l'stu}Q_{ijk}^{l'} = Q_{ij'lu}Q_{ks't}^{l}Q_{js}^{j's'}.$$

%The resulting expression contains operation $m_{4,0}$ and sends $e_{ijkstu}$ to a scalar.

%Take the scalar product of $h_{ls}$ with (11) expressing applied to $e_{ijkl}$ axiom (iii):

%$$ Q_{s'}^{qrs}Q_{ijk}^{s'} = Q_{P^2(k)t'l'}^{r}Q_{ij'}^{qt'}Q_{j}^{j'l's}.$$

%The resulting expression contains operation $m_{1,3}$ and acts on $e_{ijk}$. Making an equivalent replacement $Q_{P^2(k)t'l'}^{r} = Q_{P(t')P(l')k}^{P(r)}$, take the scalar product of it with $h^{kt}h^{ju}$:

%$$ Q_{s'}^{qrs}Q_{i}^{s'tu} = Q_{P(t')P(l')}^{P(r)t}Q_{ij'}^{qt'}Q^{j'l'su}.$$

%The resulting expression contains operation $m_{0,4}$ and acts on $e_{i}$.

%\end{proof}

All the above proves the following

\textbf{Theorem. } Any Frobenius 3-algebra with $P$ being a permutation on basis vectors has a full 3-algebra structure. 

\subsection{The case $P = id$}

Let us consider $P = id$. It is equivalent to a trivial edges labeling. Then the axioms of a Frobenius 3-algebra are simplified to the following:

(i$^*$) $(1 \otimes \overline{m})\sigma_{12}(1 \otimes \overline{m})=(\overline{m} \otimes 1)(1 \otimes \overline{m})(\overline{m} \otimes 1) : Q_{ik'}^{sq}Q_{jk}^{rk'} = Q_{i'p}^{rs}Q_{j'k}^{pq}Q_{ij}^{i'j'},$

(ii$^*$) $\overline{m}\sigma_{12}=\sigma_{12}\overline{m}: Q_{ji}^{st} = Q_{ij}^{ts},$

while a bilinear form turns out to be symmetric.

Let us fix the following notation:

1. A face is denoted by a triplet ($<$naming by vertices$>$, $<$sign of this face$>$, $<$number of a basis element corresponding to this face$>$)

2. A tetrahedron is denoted by a tetraplet consisting of its faces' triplets.

Axiom (ii$^*$) reflects the natural order of two input and two ouput faces of a tetrahedron:

1. The first input's ($ABD$) and second output's ($BDC$) common edge is considered the third vector $\overline{BD}$.

2. Choose one edge from each input face and consider them the first and the second vector so that all three vectors begin ; and constitute a right triple $(\overline{BA}, \overline{BC}, \overline{BD})$.

3. A face without this triple is considered the first output.

Thus, this tetrahedron is described by a tetraplet of triplets $((ABD, -,i), BDC, -,j), (ACD, +,s), (ABC, +,t))$ (LHS of figure 11).

This order is compatible with natural order of inputs of $m_{3,1}$ operation expressed the following way:

$$Q_{ijk}^s = Q_{jki}^s = Q_{kij}^s$$

(seen from the vertice opposite to the output face, input faces appear in the counterclockwise order).

Thus, by taking the scalar product of  $Q_{ij}^{st}$ with $h_{tk}$,  one obtains $Q_{ijk}^{s}$ denoted by tetrahedron $((ABD, -,i), BDC, -,j), (ABC, +, s), (ACD, -,k))$ (RHS of figure 11).

Again, this order of input arguments of $m_{3,1}$ is natural: looking from the vertex opposite to the output face, input facets appear in the counterclockwise order.

\begin{figure}[!ht]
\begin{center}
\includegraphics[scale=0.5]{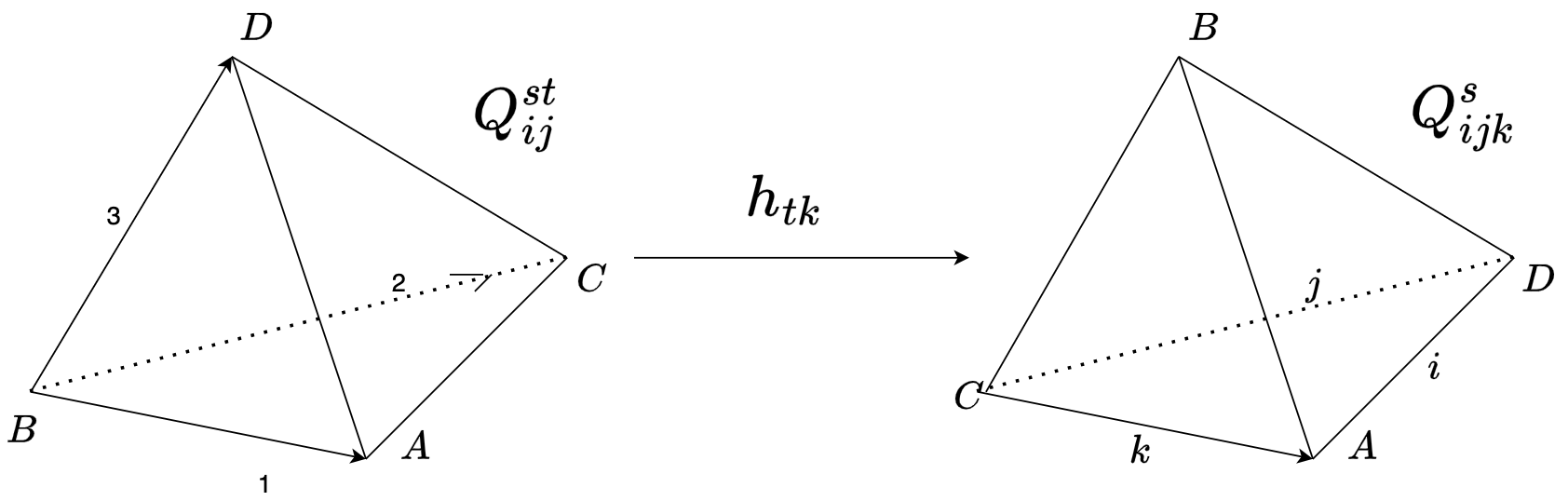}\caption{Order of faces and corresponding arguments}\label{figure11}
\end{center}
\end{figure}

%%\clearpage

2-3 Pachner move ([3]) is formalised as follows:

$$\left[
  \begin{aligned}
   ((ABC, -,i), ABD, -,j), (BDC, +,k), (ACD, +,l'))  \\
   ((ACD, -,l'), AKD, -,r), (KDC, +,s), (ACK, +,t))  \\
 \end{aligned}
  \right.$$
 
$$ \downarrow $$

$$\left[
  \begin{aligned}
   ((ABC, -,i), ABK,-,i'), (BCK, +,t'), (ACK, -,t))  \\
   ((ABD, -,j), (AKD, -,r), (BDK, +,j'), (ABK, +,i'))  \\
   ((BDK, -,j'), (BCK, -,t'), (KDC, +,s), (DBC, +,k)) \\
 \end{aligned}
  \right.$$

It provides us with a more convenient  formula rather than axiom (i$^*$) equivalent to it:

$$Q_{rl'}^{ts}Q_{ji}^{l'k} = Q_{j't'}^{sk}Q_{i'i}^{tt'}Q_{jr}^{j'i'}. \quad (15)$$

In operator form:

$$(\overline{m})_{12}(\overline{m})_{23} = (\overline{m})_{23}(\overline{m})_{13}(\overline{m})_{12} \quad (16)$$

($\overline{m}_{ij}$ acts on $i$-th and $j$-th components of a tensor product).

It is a pentagon equation ([4]).

\begin{figure}[!ht]
\begin{center}
\includegraphics[scale=0.5]{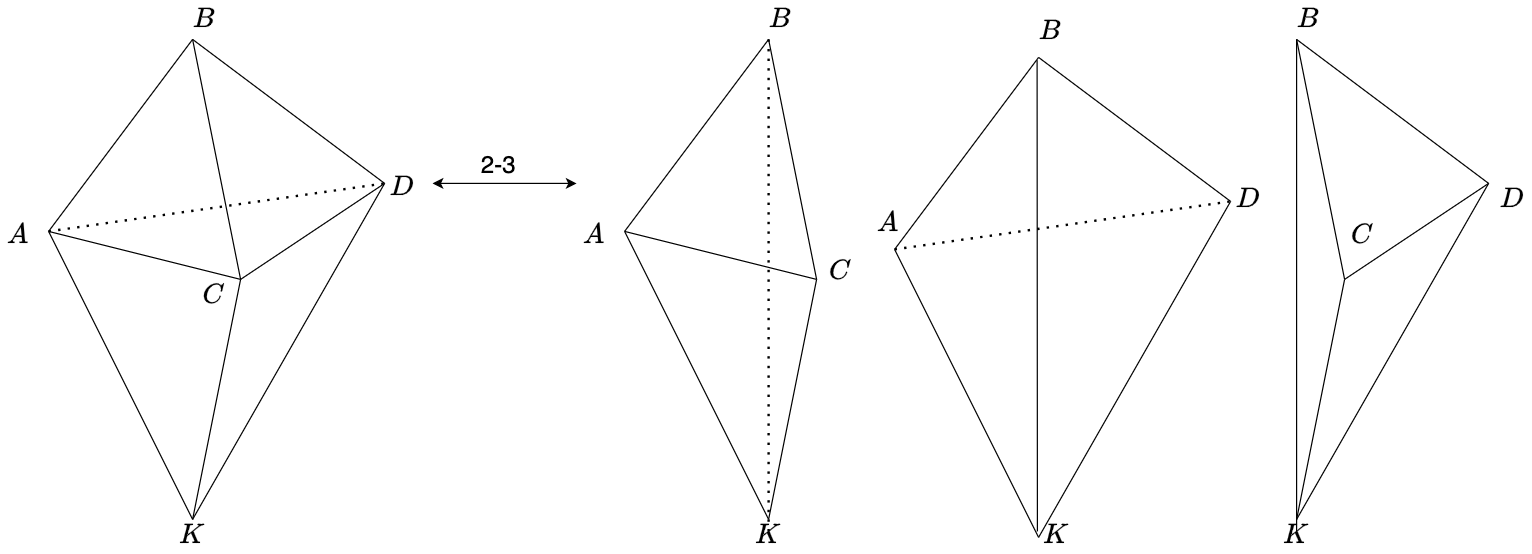}\caption{2-3 Pachner move}\label{figure12}
\end{center}
\end{figure}

%The Frobenius 3-algebra definition and Lawrence’s definition mutually harmonize conceptually—yet they are supposed to be mathematically independent. Neither specifies an explicit condition for independence under the 1-4 Pachner move. We bridge this gap by introducing such a condition. It is expressed the following way:

Lawrence's definition and Frobenius 3-algebra definition are similar but not equivalent. Additionally, none of them provides an explicit condition of independence on the 1-4 Pachner move. We bridge this gap by introducing such a condition. It is expressed the following way:

$$((ABC, -,i), (BDC, -,j), (BDC, +,k), (ABC, +,l)) $$

$$\downarrow$$

$$\left[
  \begin{aligned}
   ((ABC, -,i), AOC, +,l'), (BOC, +,k'), (ABO, -,j'))  \\
   ((ABD, -,j), (BDC, -,r'), (BDC, +,s'), (ABC, +,j'))  \\
   ((ABC, +,k), (BDC, -,s'), (BDC, +,p'), (ABC, -,k')) \\
   ((ABC, +,l), (BDC, -,l'), (BDC, -,p'), (ABC, +,r'))  \\
 \end{aligned}
  \right.$$

In coordinate form:

$$Q_{ij}^{kl} = Q_{ij'}^{k'l'}Q_{jr'}^{s'j'}Q_{k's'}^{kp'}Q_{l'p'}^{lr'} (17).$$

\textbf{Proposition. } If independence on 2-3 Pachner move holds, independence on 1-4 Pachner move may be expressed by a cubic equation.

\begin{proof}

In the RHS of (17), there is pentagon equation:

$$Q_{r'z'}^{l'p'}Q_{ji}^{z'k} = Q_{s'k'}^{p'k}Q_{j'i}^{l'k'}Q_{jr'}^{s'j'}. \quad (18)$$

Substitute (18) into (17):

$$Q_{ij}^{kl} = Q_{r'z'}^{l'p'}Q_{ji}^{z'k}Q_{l'p'}^{lr'}.  \quad (19)$$

Hence a cubic equation is obtained. 

\end{proof}

In general, one cannot simplify any further the RHS of (19 )but for some coefficients of $m_{2,2}$ it is still possible. It is sufficient to apply contraction by two free indices:

$$Q_{rl'}^{ts}Q_{ji}^{l'k} = Q_{j't'}^{sk}Q_{i'i}^{tt'}Q_{jr}^{j'i'} \rightarrow Q_{s'l'}^{ts'}Q_{ji}^{l'k} = Q_{j't'}^{s'k}Q_{i'i}^{tt'}Q_{js'}^{j'i'}. \quad (20)$$

Substitute (20) into (19):

$$Q_{i'j}^{ki'} = Q_{s'l'}^{ks'}Q_{j'j}^{l'j'}. \quad (21)$$

Let $B$ be a matrix  $\{Q_{i'j}^{ki'}\}_{j,k}$, then (21) may be rewritten as $B = B^2$ i.e. $B$ is a projector matrix.

Thus, a new family of pentagon equation solutions is introduced: arising from projectors.

\section{Appendix}

Here I introduce some examples of a 3-algebra and calculate invariants of some lens spaces in these 3-algebras.

Let us focus on $P = id$ and $B = \begin{pmatrix}
1 & 0 \\
0 & 1 \\
\end{pmatrix}$.

Some elements of $Q_{i'j}^{ki'}$:

$$Q_{1,1}^{2,1} +  Q_{2,1}^{2,2} =0,$$

$$Q_{1,2}^{1,1} + Q_{2,2}^{1,2} = 0.$$

Let us introduce variables for the structure constants of  $m_{2,2}$:

$$Q_{1,1}^{1,2} = Q_{1,1}^{2,1} = a, \quad Q_{1,2}^{1,1} = Q_{2,1}^{1,1} = b, \quad Q_{1,2}^{1,2} = Q_{2,1}^{2,1} = c, \quad Q_{1,2}^{2,1} = Q_{2,1}^{1,2} = d, \quad Q_{1,1}^{2,2} = f, \quad Q_{2,2}^{1,1} = y,$$

then:

$$Q_{1,2}^{2,2} = Q_{2,1}^{2,2} = -Q_{1,1}^{2,1} = -a, \quad  Q_{2,2}^{2,1} = Q_{2,2}^{1,2} = -Q_{1,2}^{1,1} = -b,$$

$$Q_{1,1}^{1,1} + Q_{2,1}^{1,2} = 1 \implies Q_{1,1}^{1,1} = 1 - d ,$$

$$Q_{1,2}^{2,1} + Q_{2,2}^{2,2} = 1 \implies  Q_{2,2}^{2,2} = 1 - d.$$

\textbf{1-4 Pachner move. }  Equations (19) in variables $\{ a,b,c,d,f,y\}$ are reduced to the following equation:

$$ 4ab + 2cd + (1-d)^2 + fy = 1.   \quad (22)$$

\clearpage

\textbf{Pentagon equation. }  Substituing $\{ a,b,c,d,f,y\}$ into the pentagon equation (15), one obtains the following system:

$$\begin{cases}
  a(c - d) = 0 \\
  b(c - d) = 0 \\
  (c - 1)(c - d) = 0\\
  d(c - d) = 0 \\
  f(c - d) = 0 \\
  h(c - d) = 0 \\
   ad - bf = 0 \\
   ay - bc = 0 \\
  d^2 - fy = 0 \\
  2ab + 2d^2 - d = 0 \\
  2a^2 + 2df - f = 0 \\
  2b^2 + 2dh - y = 0 \\
 \end{cases} \quad (23)$$

Notice that (23) solves (22) itself.

Possible solutions of (23):

$$\{a = b = d = f = y = 0, c = 1\},$$

$$\{a = \pm \alpha \sqrt{\frac{d}{2\alpha} - \frac{d^2}{\alpha}} , b =  \pm \sqrt{\frac{d}{2\alpha} - \frac{d^2}{\alpha}}, c = d, d, f =  \alpha d, y = \frac{d}{\alpha} \quad |  d,\alpha \in \mathbb{C^*} \},$$

$$\{a = \mp \alpha \sqrt{\frac{d^2}{\alpha} - \frac{d}{2\alpha}}, b = \pm \sqrt{\frac{d^2}{\alpha} - \frac{d}{2\alpha}}, c = d, d, f = - \alpha d, y = - \frac{d}{\alpha} \quad |  d,\alpha \in \mathbb{C^*} \}.$$

\clearpage

3-algebras obtained:

$$ \begin{cases}
  Q_{1,1}^{1,2} = Q_{1,1}^{2,1} = \pm \alpha \sqrt{\frac{d}{2\alpha} - \frac{d^2}{\alpha}} \\
  Q_{1,2}^{1,1} = Q_{2,1}^{1,1} = \pm \sqrt{\frac{d}{2\alpha} - \frac{d^2}{\alpha}} \\
  Q_{1,2}^{1,2} = Q_{2,1}^{2,1} = Q_{1,2}^{2,1} = Q_{2,1}^{1,2} = d \\
  Q_{1,1}^{2,2} = \alpha d \\
  Q_{2,2}^{1,1} = \frac{d}{\alpha} \\ 
  Q_{1,2}^{2,2}  = Q_{2,1}^{2,2} = \mp \alpha \sqrt{ \frac{d}{2\alpha} - \frac{d^2}{\alpha} } \\
  Q_{2,2}^{2,1} = Q_{2,2}^{1,2}  = \mp \sqrt{\frac{d}{ 2\alpha} - \frac{d^2}{\alpha} } \\
  Q_{1,1}^{1,1} = 1 - d \\
  Q_{2,2}^{2,2} = 1 - d \\
 \end{cases} \quad (24)$$

A lens space $L(p,q)$ may be  split into tetrahedra via a bipyramide  over a $p$-gon $A_1...A_p$ with poles $N$ and $S$ (for conveniency, let us say $A_{p+1} = A_p$): 

1. split a $p$-gon $A_1...A_p$ into triangles using diagonals from $A_p$ to other vertices

2. construct two tetrahedra over each obtained triangle: one from $N$, another from $S$

3. identify $A_lA_{l+1}N$ with $A_{(l + q \quad mod \quad p)}A_{(l + 1 + q \quad mod \quad p)}$

(see figure 13).

\begin{figure}[!ht]
\begin{center}
\includegraphics[scale=0.6]{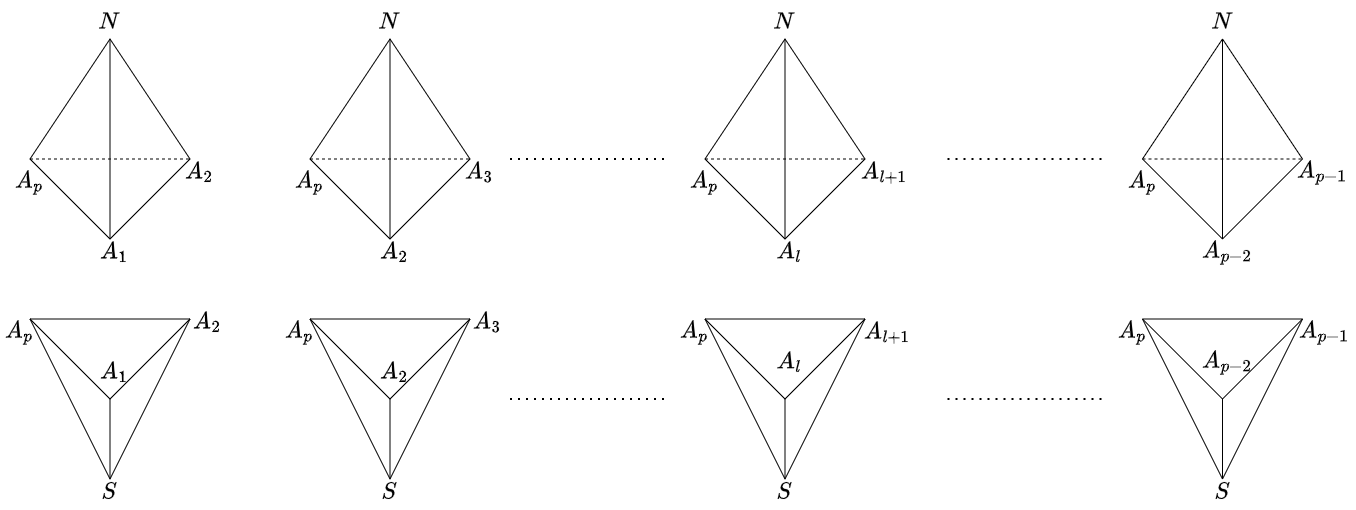}\caption{Lens space splitting}\label{figure13}
\end{center}
\end{figure}

\clearpage

Splitting $L(3,1)$:

$$\left[
  \begin{aligned}
   ((A_1A_{2}N, -,i_1), (A_1A_3N, -,j_1), (A_{2}A_3N, +,j_{2}), (A_1A_{2}A_3, +,t_1))  \\
   ((A_1A_{2}S, +,j_1), (A_1A_3S, -,j_2), (A_{2}A_3S, +,t_1), (A_1A_{2}A_3, - ,t_1))  \\
 \end{aligned},
  \right.$$

Hence $L(3,1)$ is assigned the value $Q_{i_1j_1}^{j_2t_1}Q_{j_2t_1}^{j_1i_1}$.

Solely relying on (23) , calculate it:

 $$Q_{i_1j_1}^{j_2t_1}Q_{j_2t_1}^{j_1i_1} = 8ab + 6d^2 - 4d + 2fy +2 = 2.$$

Splitting $L(4,1)$:

$$\left[
  \begin{aligned}
   ((A_1A_{2}N, -,i_1), (A_1A_4N, -,j_1), (A_{2}A_4N, +,j_{2}), (A_1A_{2}A_4, +,t_1))  \\
   ((A_1A_{2}S, +,j_1), (A_1A_4S, -,j_3), (A_{2}A_4S, +,s_1), (A_1A_{2}A_4, - ,t_1))  \\
    ((A_2A_{3}N, -,i_2), (A_2A_4N, -,j_2), (A_{3}A_4N, +,j_{3}), (A_2A_{3}A_4, +,t_2))  \\
   ((A_2A_{3}S, +,i_1), (A_2A_4S, -,s_1), (A_{3}A_4S, +,i_2), (A_2A_{3}A_4, - ,t_2))  \\
 \end{aligned},
 \right.$$

Hence $L(4,1)$ is assigned a value $Q_{i_1j_1}^{j_2t_1}Q_{i_2j_2}^{j_3t_2} Q_{j_3t_1}^{j_1s_1}Q_{s_1t_2}^{i_1i_2}$.

Solely relying on (23) , calculate it:

 $$Q_{i_1j_1}^{j_2t_1}Q_{i_2j_2}^{j_3t_2} Q_{j_3t_1}^{j_1s_1}Q_{s_1t_2}^{i_1i_2} = 4a^2h + 2a^2 - 8abd +12ab + 4b^2f + 2b^2  - 4d^3 + 12d^2 + 4dfy + 2df  + 2dy - $$
 
 $$ -8d + 2   = 2 - 2d + f + y. \quad (25)$$
 
 Calculate (25) within the family of 3-algebras (24):
 
 $$Q_{i_1j_1}^{j_2t_1}Q_{i_2j_2}^{j_3t_2} Q_{j_3t_1}^{j_1s_1}Q_{s_1t_2}^{i_1i_2} =  2 - 2d + f + y = 2 - 2d + \alpha d + \frac{d}{\alpha} = $$
\\
$$= 2 + (\sqrt{\alpha} - \frac{1}{\sqrt{\alpha}})^2d$$

To distinguish between homotopically non-equivalent spaces $L(3,1)$ and $L(4,1)$ ($\pi_1(L(3,1)) = \mathbb{Z}_3$, $\pi_1(L(4,1)) = \mathbb{Z}_4$), it suffices to set $\alpha \neq \pm 1$.

\section{Further questions}

In this paper, I heavily rely on invertibilty of a bilinear form. What can be if one takes a degenerate one? 

Are there 3-algebras which have no natural full 3-algebra structure? What can one say about them?

How to deal with an 3-algebra whose $P$ cannot be represented as a permutation on basis vectors?

How to generalise Frobenius compatibility to $n$-algebras?

What is the interplay between lens spaces and 3-algebras?

What can be said about pentagon equation solutions arising from projectors?

% ===== Bibliography =====

\end{document}